\documentclass[]{iopart}
\usepackage{iopams}
\usepackage{setstack}
\usepackage[]{graphicx}
\usepackage{color}

\usepackage[latin1]{inputenc}
\usepackage[english]{babel}
\usepackage{dsfont}

\newcommand{\by}{\bi{y}}
\newcommand{\ba}{\bi{a}}
\newcommand{\rT}{{\rm T}}
\newcommand{\ha}{\hat{\bi{a}}}
\newcommand{\hy}{\hat{\bi{y}}}
\newcommand{\he}{\hat{\bepsilon}}
\newcommand{\Prob}{{\rm Prob}}

\begin{document}

\title[]{Bayesian estimate of the degree of a polynomial given a noisy data sample}
\author{Giovanni Mana, Paolo Alberto Giuliano Albo, and Simona Lago}
\address{INRIM -- Istituto Nazionale di Ricerca Metrologica, str.\ delle Cacce 91, 10135 Torino, Italy}

\begin{abstract}
A widely used method to create a continuous representation of a discrete data-set is regression analysis. When the regression model is not based on a mathematical description of the physics underlying the data, heuristic techniques play a crucial role and the model choice can have a significant impact on the result. In this paper, the problem of identifying the most appropriate model is formulated and solved in terms of Bayesian selection. Besides, probability calculus is the best way to choose among different alternatives. The results obtained are applied to the case of both univariate and bivariate polynomials used as trial solutions of systems of thermodynamic partial differential equations.
\end{abstract}

\submitto{Metrologia}
\pacs{07.05.Kf, 02.50.Cw, 06.20.Dk, 02.50.Tt}

% 06.20.Dk Measurement and error theory
% 02.50.Cw Probability theory
% 07.05.Kf Data analysis: algorithms and implementation; data management
% 02.50.Tt Inference methods

\ead{g.mana@inrim.it}

\section{Introduction}
A problem in regression analysis is to determine how many basis functions to include in the regression model, for instance, when determining the calibration curve that best fits the data \cite{Massa:2013}. Any set of basis functions can be considered; when they are polynomials, the problem is determining the degree of the regression. A maximum likelihood approach, which leads to the highest possible number of the basis functions, cannot be the right choice. This problem has been considered by many authors in different statistical settings and their investigations led to a number of proposed solutions \cite{Anderson:1961,Schartz:1978,Shao:1996,Philips:1998,Guttman:2005}. An original and undeservedly neglected one is hidden in a tutorial paper on Bayesian reasoning by Gull \cite{Gull:1988}, where the basic idea is to calculate and to compare each model probability, given the data.

In order to bring this result to the metrologist's attention, we reassess the Gull work and make clear its usefulness in selecting among linear models. In addition, by slightly changing the model parametrisation, we obtain an exact analytical solution and demonstrate that, in suitable limit cases, it reduces to the Gull's approximate one.

Here obtained results may have an impact when the polynomial coefficients are used for solving systems of partial differential equations \cite{Lago:2008,Lago:2013}. In this case, different choices of the polynomial degree lead to different sets of coefficients and, consequently, to different solutions. The availability of a rigorous criterion based on the probability calculus allows any arbitrary choice -- in general, driven only by the residuals analysis -- to be avoided. To illustrate the concepts here described, it is shown how to determine the set of basis functions that best fits the measured values of the speed of sound in acetone, as a function of the temperature and pressure.

\section{Problem statement}
We want to represent the $\by=[y_1, y_2, ...\, y_N]^\rT$ measurement results by the linear model
\begin{equation}\label{model}
 \by = W(l) \ba + \bepsilon ,
\end{equation}
where $\bepsilon=[\epsilon_1, \epsilon_2, ...\, \epsilon_N]^\rT$ are additive uncorrelated Gaussian errors having unknown variance $\sigma^2$ and zero mean, $\ba=[a_0, a_1, ...\, a_{l-1}]^\rT$ are $l$ model parameters, $W(l)$ is a $N\times l$ matrix explaining the data, $W_{nm}=w_m(x_n)$, and $\{w_0(x), w_1(x), ...\, w_{l-1}(x)\}$ is a set of $l$ basis functions. The basis functions may be polynomials, for instance, $w_m(x)=x^m$, but, in general, they are any set of linearly independent functions. The problem is to find the set of basis functions most supported by the data; when they are polynomials, this corresponds to find the optimal degree of the regression.

The interpretative model of the data is summarised by the matrix $W(l)$; therefore, the problem is equivalent to find -- within a set of matrices explaining the data -- the one most supported by the data. Since it explicitly appears in the final formulae and for the sake of notational simplicity, we label the $W$ matrices by the number $l$ of the model free-parameters. However, we can compare also models having the same number of parameters but different basis functions.

\section{Bayesian inferences}
According to the Bayes theorem -- by assigning the same probability to all the models -- the probability of the $l$-th model to explain the data is proportional to the probability of the observed data given $W(l)$, no matter what the values of the model parameters may be. In turn, it is the normalising factor of the likelihood of the model parameters times the probability distribution synthesising the information available about the parameter values before the measurement results are available.

To steer the calculation, we must first determine the post-data probability density, $P(\ba,\sigma|\by,l)$, of the parameters of each model (which parameters include the unknown variance $\sigma^2$ of the data) given the $y_n$ data and the data-explaining matrix $W(l)$. The post-data probability density is found via the product rule of probabilities \cite{Sivia,Jaynes},
\begin{equation}\label{product}
 P(\ba,\sigma|\by,l) Z(\by|l) = N_N(\by|\sigma,\ba,l) \pi(\sigma,\ba|l) ,
\end{equation}
where the $N$-dimensional Gaussian function
\begin{equation}\label{lik}
 N_N(\by|\sigma,\ba,l) = \frac{1}{\sqrt{(2\pi)^N}\,\sigma^N} \exp \left( -\frac{|\by-W\ba|^2}{2\sigma^2} \right)
\end{equation}
is the likelihood of the $\ba$ and $\sigma$ parameters, $\pi(\sigma,\ba|l)$ is their pre-data probability density, the sought normalisation factor of $N_N(\by|\sigma,\ba,l) \pi(\sigma,\ba|l)$,
\begin{equation}\label{evidence}
 Z(\by|l) = \int_\mathit{\Gamma} N_N(\by|\sigma,\ba,l) \pi(\sigma,\ba|l)\, \rmd \sigma \rmd \ba ,
\end{equation}
is named model evidence, and the integration is carried out over the hypervolume $\mathit{\Gamma}$ associated to the possible $\ba$ and $\sigma$ values.

Next, by observing that $Z(\by|l)$ is also the probability density of the data given $W(l)$ -- whatever the values of $\ba$ and $\sigma$ may be -- we get the post-data model-probability, $\Prob(l\,|\by)$, by applying again the product rule of probabilities to the $\{l, \by\}$ pair. Hence \cite{McKay},
\begin{equation}\label{degree}
 \Prob(l\,|\by) A(\by) = Z(\by|l) ,
\end{equation}
where, prior the data are at hand, we assigned the same probability to each model and
\begin{equation}\label{norm}
 A(\by) = \sum_{l=1}^L Z(\by|l) ,
\end{equation}
where $L$ is the number of models to be compared, is the normalisation factor of $Z(\by|l)$. Therefore, to solve the stated problem, the calculation of the evidence (\ref{evidence}) is central.

\section{Pre-data distribution}
To set the pre-data distribution of $\ba$ and $\sigma$, we assume that they are independent. Hence, $\pi(\sigma,\ba|l)=\pi_\sigma(\sigma) \pi_a(\ba|l)$. As regards $\sigma$, we use the improper Jeffreys prior \cite{Jaynes:1968}
\begin{equation}\label{Jeffreys}
 \pi_\sigma(\sigma) = 1/\sigma ,
\end{equation}
which is invariant for a change of the measurement unit of the data.

As regards the $\ba$ parameters, let the mean of $\by$, whatever the $\ba$ values may be, null. The relevant average is carried out over the joint distribution $N_N(\by|\sigma,\ba,l) \pi(\sigma,\ba|l)$, not over the sampling distribution of the data $N_N(\by|\sigma,\ba,l)$, where the $\ba$ values are fixed. Consequently, since $\by = W \ba + \bepsilon$ and $\bepsilon$ are zero-mean errors, also the pre-data mean of the $\ba$ parameters is zero.

In addition, let $\beta^2\mathds{1}$ be the pre-data covariance of $\by$, whatever the $\ba$ values may be and where $\mathds{1}$ is the unit matrix. Also in this case, the stated covariance is relevant to the joint distribution $N_N(\by|\sigma,\ba,l) \pi(\sigma,\ba|l)$. Hence, by observing that $W C_{aa} W^\rT + \sigma^2 \mathds{1} = \beta^2\mathds{1}$, because of (\ref{model}), the pre-data covariance $C_{aa}$ of the model parameters is
\begin{equation}\label{Caa}
 C_{aa} = (\beta^2-\sigma^2) (W^\rT W)^{-1} ,
\end{equation}
where $\beta^2 > \sigma^2$.

Eventually, since the prior distribution of $\ba$ is constrained by $\langle \ba \rangle = 0$ and (\ref{Caa}), where the angle bracket indicate the mean, the principle of maximum entropy fixes the sought pre-data distribution to the $l$-dimensional Gaussian distribution \cite{Jaynes}
\begin{equation}\label{prior:a}
 \pi_a(\ba|\beta,l) = N_l(\ba|0,C_{aa})
\end{equation}
having zero mean and $C_{aa}$ covariance. Actually, the $\beta$ value is unknown. Therefore, we should eliminate it by marginalisation,
\begin{equation}\label{prior:c}
 \pi_a(\ba|l) = \int_\sigma^\infty \pi_a(\ba,\beta|l)\, \rmd \beta ,
\end{equation}
where
\begin{equation}\fl\label{prior:b}
 \pi_a(\ba,\beta|l) = \pi_a(\ba|\beta,l)\pi_\beta(\beta) =
 \frac{1}{\beta} \sqrt{ \frac{\det(W^\rT W)}{(2\pi)^l\,(\beta^2-\sigma^2)^l} }\,
 \exp \left[ -\frac{\ba^\rT W^\rT W \ba}{2(\beta^2-\sigma^2)} \right] ,
\end{equation}
and $\pi_\beta(\beta)=1/\beta$ is a Jeffreys prior. However, since the integration in (\ref{prior:c}) cannot be done analytically, we add $\beta$ to the model parameters and delay the marginalisation over $\beta$ as much as possible. Hence, by putting (\ref{Jeffreys}) and (\ref{prior:b}) together, we obtain the pre-data distribution of the full set of model parameters,
\begin{equation}\label{prior:d}
 \pi(\sigma,\beta,\ba|l) = \frac{1}{\beta\sigma}
 \sqrt{ \frac{\det(W^\rT W)}{(2\pi)^l\,(\beta^2-\sigma^2)^l} }\,
 \exp \left[ -\frac{\ba^\rT W^\rT W \ba}{2(\beta^2-\sigma^2)} \right] .
\end{equation}

A comment is necessary. In general, the use of improper priors -- like $\pi_\sigma(\sigma) = 1/\sigma$ and $\pi_\beta(\beta) = 1/\beta$ -- should be avoided, because, in such a case, the model evidence (\ref{evidence}) is defined only up to unknown scale factors. However, since in this case the same factor is included in all the evidences, this does not jeopardise the model comparison.

\section{Model selection}
\subsection{Evidence calculation}
By combining (\ref{lik}) and (\ref{prior:d}), the evidence of $W(l)$ is
\begin{eqnarray}\nonumber
 \fl Z(\by|l) = \int_0^{+\infty}\rmd\sigma \frac{1}{\sigma^{N+1}}
 \int_\sigma^{+\infty} \rmd\beta \frac{1}{\beta} \sqrt{\frac{\det(W^\rT W)}{(2\pi)^{N+l}\,(\beta^2-\sigma^2)^l}} \\
 \lo\times \int_{-\infty}^{+\infty} \exp \left[ -\frac{|\by-W\ba|^2}{2\sigma^2} \right]
           \exp \left[ -\frac{\ba^\rT W^\rT W \ba}{2(\beta^2-\sigma^2)} \right] \rmd\ba .
\end{eqnarray}

Before carrying out the integration, we observe that
\begin{equation}\label{formula}
 |\by-W\ba|^2 = (\ba-\ha)^\rT W^\rT W (\ba-\ha) + \by^\rT (\by - \hy) ,
\end{equation}
where $\ha=(W^\rT W)^{-1} W^\rT \by$ is the least squares estimate of $\ba$ and $\hy = W\ha$ is the measurand estimate. Hence, the first integration is
\begin{eqnarray}\nonumber\fl
 \exp \left[ \frac{-\by^\rT \he}{2\sigma^2} \right]
 \int_{-\infty}^{+\infty} \exp \left[ \frac{-(\ba-\ha)^\rT W^\rT W (\ba-\ha)}{2\sigma^2} \right]
 \exp \left[ \frac{-\ba^\rT W^\rT W \ba}{2(\beta^2-\sigma^2)} \right]\, \rmd \ba = \\
 \left( \frac{\sigma}{\beta} \right)^l \sqrt{ \frac{(2\pi)^l\,(\beta^2-\sigma^2)^l}{\det(W^\rT W)} }\,
 \exp \left[ \frac{-\by^\rT \he}{2\sigma^2} \right]
 \exp \left[ \frac{-|\hy|^2}{2\beta^2} \right] ,
\end{eqnarray}
where $\he=\by-\hy$ are the residuals and $|\hy|^2=\hy^\rT \hy=\ha^\rT W^\rT W \ha$.  It must be noted that $\by^\rT \he > 0$ because $|\by-W\ba|^2 > 0$ no matter what the $\ba$ value may be. Consequently, the right-hand side of (\ref{formula}) is greater than zero also when $\ba=\ha$. Hence, $\by^\rT \he = \by^\rT(\by-\hy) > 0$.

The next integration,
\begin{eqnarray}\nonumber
 \frac{1}{\sigma^{N+1-l} \sqrt{(2\pi)^N}} \exp \left[ \frac{-\by^\rT \he}{2\sigma^2} \right]
 \int_\sigma^{+\infty} \frac{1}{\beta^{l+1}}
 \exp \left[ \frac{-|\hy|^2}{2\beta^2} \right] \rmd\beta = \\
 \frac{\sqrt{2^{l-2}} \left\{ \Gamma(l/2) - \Gamma\big[p/2,|\hy|^2/(2\sigma^2)\big] \right\} }
 {\sigma^{N+1-l}\sqrt{(2\pi)^N} \, |\hy|^l}
 \exp \left[ \frac{-\by^\rT \he}{2\sigma^2} \right] ,
\end{eqnarray}
where $\Gamma(z)$ is the Euler gamma function, eliminates $\beta$.

Eventually, provided $N>l$, the evidence is
\begin{eqnarray}\nonumber\fl
 Z(\by|l) = \frac{\sqrt{2^{l-2}}}{\sqrt{(2\pi)^N} \,|\hy|^l}
 \int_0^{+\infty} \frac{\Gamma(l/2) - \Gamma\big[p/2,|\hy|^2/(2\sigma^2)\big]}{\sigma^{N+1-l}}
 \exp \left[ \frac{-\by^\rT \he}{2\sigma^2} \right] \rmd\sigma \\ \fl
 = \frac{\Gamma\left(\frac{N-l}{2}\right)}{4\sqrt{\pi^N}} \left[
 \frac{\Gamma\left(\frac{l}{2}\right)}{|\hy|^l (\by^\rT \he)^{(N-l)/2}}
 - \frac{\Gamma\left(\frac{l}{2}\right)\, _2{\tilde F}_1\left( \frac{l}{2}, \frac{N-l}{2}; \frac{N+2-l}{2};
 \frac{\by^\rT \he}{|\hy|^2} \right)}{|\hy|^N} \right] ,
\end{eqnarray}
where $_2{\tilde F}_1(a,b;c;z)$ is the regularised hypergeometric function.

\subsection{Model probability}
By assuming that, prior the data are available, each $W(l)$ has the same probability, according to (\ref{degree}) and (\ref{norm}), the $l$-model probability is proportional to the $l$-model evidence; that is, 
\begin{eqnarray}\nonumber
 \Prob(l|\by) \propto Z(\by|l) \propto \\ \fl
 \frac{\Gamma\left(\frac{N-l}{2}\right) \Gamma\left(\frac{l}{2}\right)}{|\hy|^l (\by^\rT \he)^{(N-l)/2}}
 - \frac{\Gamma\left(\frac{N-l}{2}\right) \Gamma\left(\frac{N}{2}\right)\,
 _2{\tilde F}_1\left( \frac{l}{2}, \frac{N-l}{2}; \frac{N+2-l}{2};
 \frac{\by^\rT \he}{|\hy|^2} \right)} {|\hy|^N} . \label{evidence:p}
\end{eqnarray}
It is worth noting that, since $Z(\by|l)$ is the marginal probability density of the data given $W(l)$ (no matter what the values of the model parameters may be), the dimensions of $\Prob(l|\by)$ are the same of $|\hy|^{-N}$.

If $W(l_0)$ explains the data exactly, then $\he(l_0)=0$ and $\Prob(l_0|\by)=1$, as expected. Furthermore, $\Prob(l|\by)$ is independent of the $\by$ scale. In fact, when $\by \rightarrow \lambda\by$, the evidence of $l$ transforms as $Z(\by|l) \rightarrow \lambda^{-N} Z(\by|l)$, which leaves $\Prob(l|\by)$ unchanged. In addition, $Z(\by|l)$ depends only on $\by$ and $\hy$; therefore, it is independent of the choice of the sampling points $x_n$. Eventually, $\Prob(l|\by) $ is not invariant for translation of the origin of the $y$-axis; this is a consequence of the $\langle \by \rangle=0$ assumption, which is embedded into the pre-data distribution of the $\ba$ coefficients.

\subsubsection{Asymptotic behaviours.}
In the case when $N-l \gg 2$, we can use the approximations $(N+2-l)/2 \approx (N-l)/2$ and
\begin{equation}
 _2{\tilde F}_1(N/2,(N-l)/2;(N-l)/2;z) \approx \frac{(1-z)^{-N/2}}{\Gamma[(N-l)/2]} .
\end{equation}
In addition, for a large data sample, since $\hy^\rT \he /\chi_y^2 \approx 0$, where $\chi_y^2 = |\hy|^2$ is the sum of the squared data, and $\by=\hy+\he$, it follows that $\by^\rT \he /\chi_y^2 \approx |\he|^2 /\chi_y^2 = \chi_\epsilon^2 /\chi_y^2$, where $\chi_\epsilon^2 =|\he|^2$ indicates the sum of the squared residuals. Therefore, apart from the $\hy^{-N} \approx$ const.\ factor that we omit, we can rewrite (\ref{evidence:p}) as
\begin{equation}\label{evidence:p1}
 \Prob(l|\by) \approx  \frac{\Gamma[(N-l)/2]\Gamma(l/2)}{(\chi_\epsilon/\chi_y)^{N-l}}
 - \frac{\Gamma(N/2)} {(1 + \chi_\epsilon^2/\chi_y^2)^{N/2}} .
\end{equation}

Eventually, if  $\chi_\epsilon^2/\chi_y^2 \ll 1$ -- which means good data and good models -- and $N\gg l$, (\ref{evidence:p1}) simplifies further as
\begin{equation}\label{evidence:p2}
 \Prob(l|\by) \approx \frac{\Gamma[(N-l)/2]\Gamma(l/2)}{(\chi_\epsilon/\chi_y)^{N-l}} .
\end{equation}
This equation brings into light that, among the models having the same number of free parameters, the most supported by the data is that whose associated sum of the squared residuals is minimum. In additions, it shows that the optimal model minimises the residuals by keeping at the same time $l$ as small as possible, in order to maximise $N-l$.

As a last step we write (\ref{evidence:p2}) as
\begin{equation}\fl
 \ln[\Prob(l|\by)] \approx \ln\left[\Gamma\left(\frac{N-l}{2}\right)\right] + \ln\left[\Gamma\left(\frac{l}{2}\right)\right] + (N-l)\ln(\chi_y/\chi_\epsilon) ,
\end{equation}
which is the approximation given in \cite{Gull:1988}.

Now, we study the case when the data are samples of a polynomial having degree $q$ and the $\bepsilon$ variance tends to zero. In this case, provided $l \ge q+1$, the residuals are independent of the degree of the fitting polynomial, $\chi_y/\chi_\epsilon \rightarrow \infty$, and
\begin{equation}
 \ln[\Prob(l|\by)] \approx (N-l)\ln(\chi_y/\chi_\epsilon) ,
\end{equation}
which shows that the evidence is maximum when $l$ is minimum. Therefore, the degree most supported by the data is $q$, as expected.

If $\chi_\epsilon/\chi_y \rightarrow \infty$ -- which means bad data or bad models -- and $N\gg l$, (\ref{evidence:p1}) simplifies as
\begin{equation}
 \Prob(l|\by) \approx \Gamma[(N-l)/2]\Gamma(l/2) ,
\end{equation}
which indicates that the optimal data model has only one degree of freedom, that is, $y_n=a_0+\epsilon_n$.

\begin{figure}
\centering
\includegraphics[width=63mm]{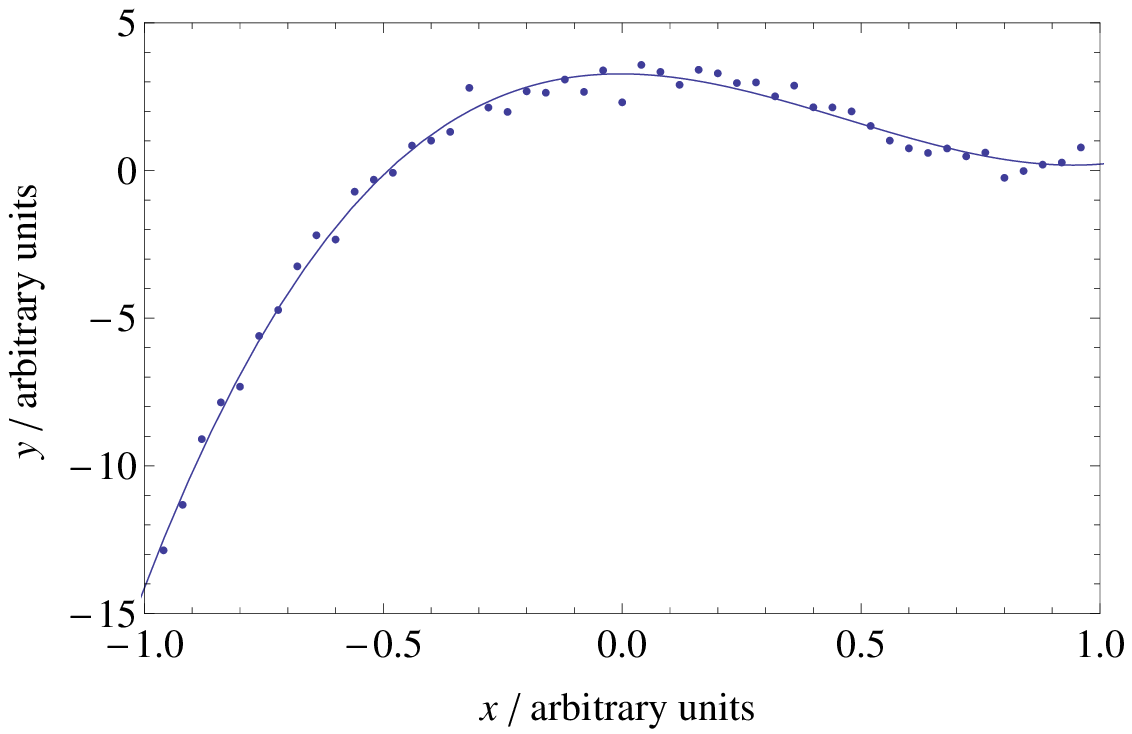}
\includegraphics[width=63mm]{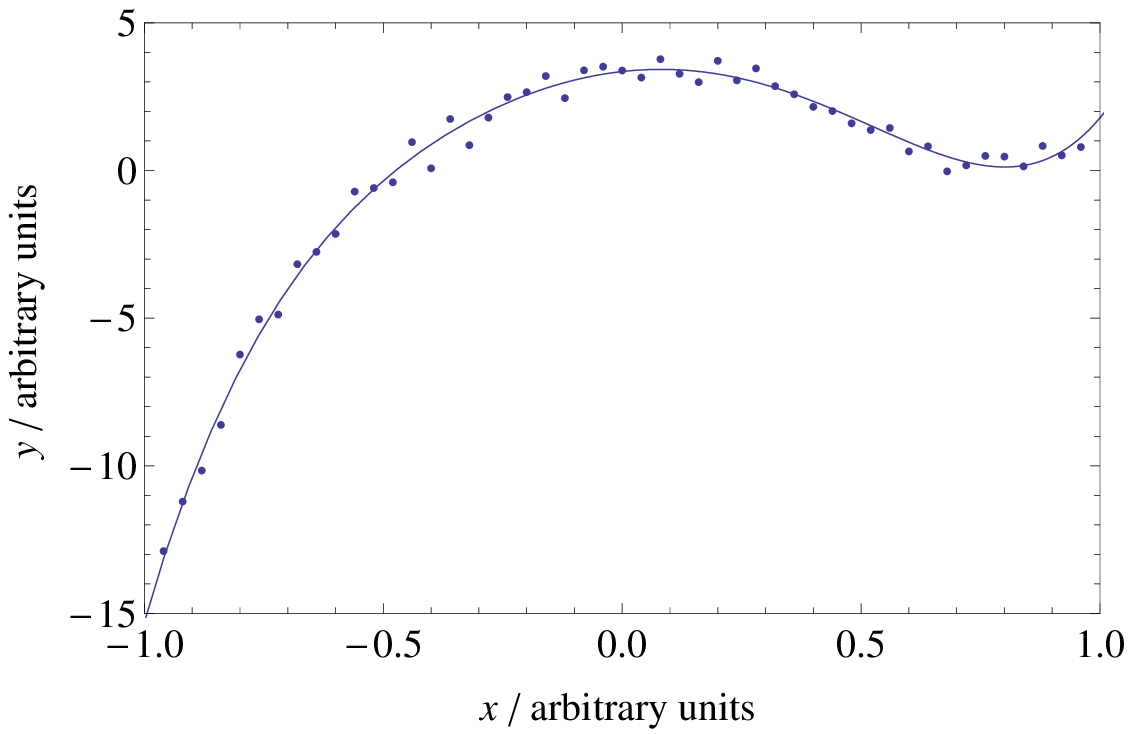}
\caption{Simulated noisy data-samples of the polynomial (\ref{poly}); solid line are the most probable polynomials explaining these specific data sets (left, a third degree polynomial; right, a fifth degree polynomial).} \label{data}
\end{figure}

\section{Numerical example}
The figure \ref{data} shows two independent sets of $N=50$ simulated data each, uniformly sampled in the $[x_{\rm min}=-1, x_{\rm max}=1]$ interval from the fifth degree polynomial
\begin{equation}\label{poly}
 y = - 1 x - 10 x^2 + 2 x^3 + 5 x^5 .
\end{equation}
The outputs of a Gaussian random-number generator having zero mean and 0.4 standard deviation were added to the data. In order to fulfill the $\langle \by \rangle = 0$ requirement, the data have been pre-processed to remove the arithmetic mean. The $\sigma=0.4$ value of the $\bepsilon$ standard deviation was chosen intermediate between the good and bad data limit cases.

\begin{figure}
\centering
\includegraphics[width=63mm]{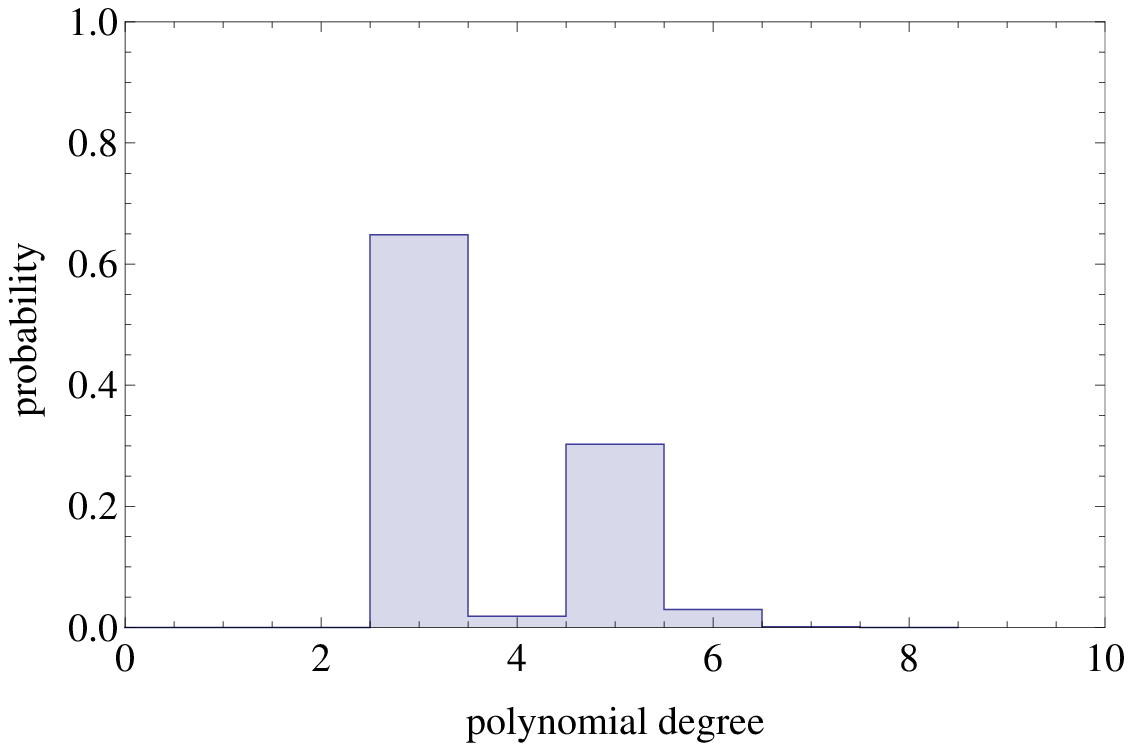}
\includegraphics[width=63mm]{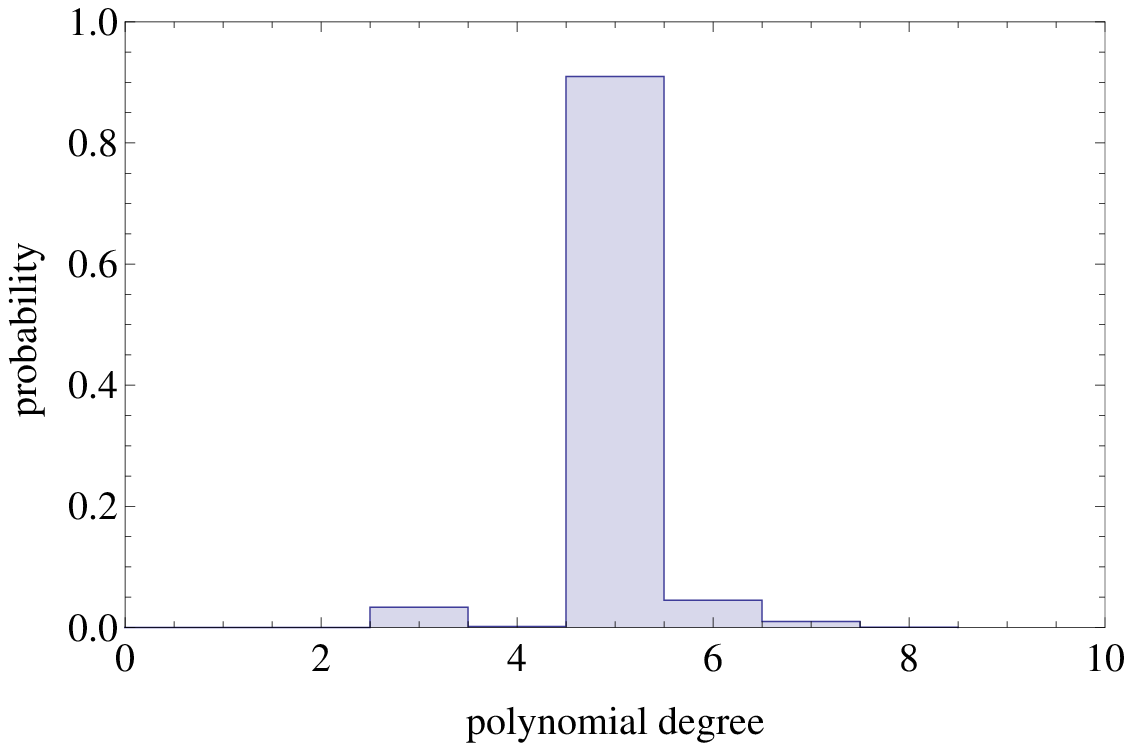}
\caption{Probabilities of the degree of the polynomials explaining the data sets shown in Fig. \ref{data}.} \label{degree-fig}
\end{figure}

To explain the data, a set of ten polynomials -- with degree from zero to nine -- have been considered, each polynomial has been fitted to the data, and both the error and data estimates -- $\he$ and $\hy$, respectively -- have been calculated. Eventually, the evidence of each polynomial has been found by application of (\ref{evidence:p}) as well as, after normalisation, each polynomial probability to explain the data has been calculated. The results for each of the data sets shown in Fig.\ \ref{data} are shown in Fig.\ \ref{degree-fig}.

With the $\sigma=0.4$ choice, the degree of the polynomials that best explain the data have been always found equal to three or five, depending on the specific data set. The figures \ref{data} and \ref{degree-fig} show these alternating cases. It is worth noting that the probability to explain the data of a polynomial of fourth degree -- whose basis function is missing in (\ref{poly}) -- is very low. As expected, by reducing the noise level, the most likely degree stabilises on five, whereas, by increasing the noise level, it stabilises on three.

\section{Speed of sound in acetone}
In \cite{Lago:2008} it was shown how to solve the thermodynamic differential equations
\numparts\begin{eqnarray}
 \bigg( \frac{\partial \rho}{\partial p} \bigg)_T &= &\frac{T}{c_p}
 \frac{1}{\rho^2} \left( \frac{\partial \rho}{\partial T} \right)_p^2 + \frac{1}{w^2}
 \label{pde1} \label{eq1} \\
 \left( \frac{\partial c_p}{\partial p} \right)_T &= &- \frac{T}{\rho}
 \left[\frac{2}{\rho^2} \left( \frac{\partial \rho}{\partial T}
 \right)_p^2 - \frac{1}{\rho}\left(\frac{\partial^2 \rho}{\partial T^2}\right)_p \right]
\label{pde2}
\end{eqnarray}\endnumparts
relating density $\rho(p,T)$, heat capacity $c_p(p,T)$, and speed of sound $w(p,T)$, as a function of the temperature $T$ and pressure $p$. These equations can be solved if initial conditions $\rho(p_0,T)$ and $c_p(p_0,T)$ are given at a the reference pressure, $p_0$, and the speed of sound is known on the entire range of pressures and temperatures of interest.

When a numerical integration of (\ref{eq1}-$b$) is carried out, the heat capacity shows often diverging values at the extremes of the temperature range. Approaching the integration problem by using local polynomial representations of the thermodynamic quantities eliminates the divergence and allows the uncertainty of the results to be estimated. Hence, by using the trial solutions
\numparts\begin{eqnarray}
 \rho(p,T) &= &\sum_{i,j} a_{i j} (p-p_0)^i(T-T_0)^j , \\
 c_p(p,T)  &= &\sum_{i,j} b_{i j} (p-p_0)^i(T-T_0)^j , \\
 w(p,T)    &= &\sum_{i,j} c_{i j} (p-p_0)^i(T-T_0)^j ,
\end{eqnarray}\numparts
and -- once the degrees of the polynomials have been fixed -- the unknown coefficients $a_{i j}$, $b_{i j}$, and $c_{i j}$ are obtained as described in \cite{Lago:2008}. Briefly, the best polynomial approximations of the initial conditions and speed of sound are used to determine $a_{0 j}$, $b_{0 j}$, and $c_{i j }$; subsequently, the remaining coefficients $a_{i j}$ and $b_{i j}$ are calculated by means of the equations (\ref{pde1}) and (\ref{pde2}).

\begin{figure}
\centering
\includegraphics[width=65mm]{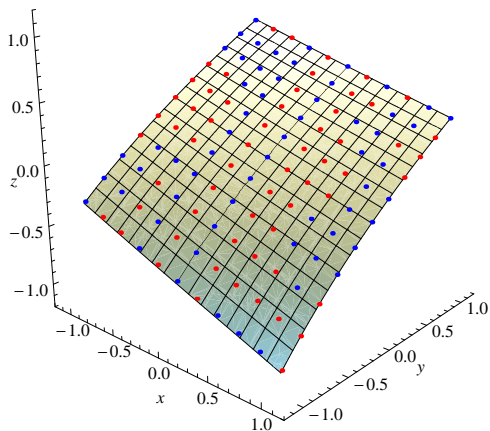}
\caption{Measured values of the speed of sound in acetone. The data have been scaled in $[-1,+1]$ intervals; $x$ is the temperature, $y$ is the pressure, and $z$ is the speed of sound. The polynomial model most likely supported by the data is also shown. Red dots are the data higher than what predicted by the model, blue dots are those lower.} \label{speed}
\end{figure}

As an application example, we show how to determine the optimal polynomial when smoothing the measured values of the speed of sound in acetone as a function of temperature and pressure \cite{Lago:2008,Lago:2013}; a set of measurement results is shown in Fig.\ \ref{speed}. For the sake of numerical convenience, the temperature, pressure, and speed have been scaled in $[-1,+1]$ intervals according to $x=(T-T_0)/\Delta T$, $y=(p-p_0)/\Delta p$, and $z=(w-w_0)/\Delta w$, where the offsets and scale factors have been suitably chosen.

\begin{figure}
\centering
\includegraphics[width=65mm]{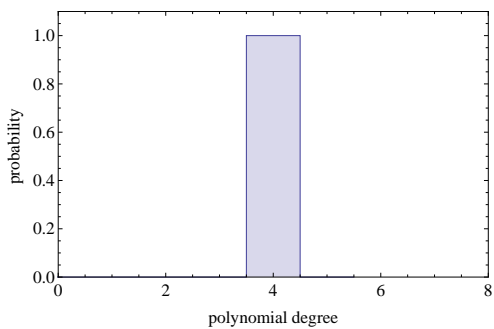}
\caption{Probability of the polynomials having degree from zero to six to explain the data in Fig. \ref{speed}.} \label{model-1}
\end{figure}

As shown in Fig.\ \ref{model-1}, the regression analyses using the seven basis-function sets $\{1, x, y, ...\, x^r y^s, ...\}_q$, where $0 \le r+s \le q$ and $q=0, 1, ...\, 6$, indicate that the $q=4$ set is the only one supported by the data. This sharp selection is due to the fast increase of the number $l$ of basis functions as the polynomial degree increases. For instance, if $q=3$, then $l=10$; if $q=4$, then $l=15$; if $q=5$, then $l=21$.

\begin{figure}
\centering
\includegraphics[width=65mm]{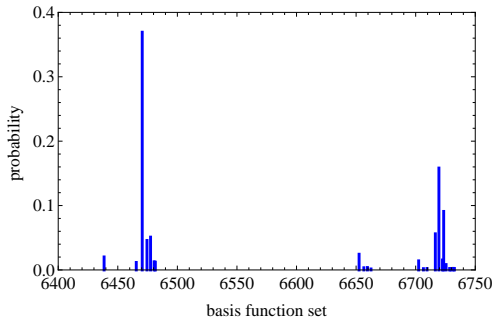}
\caption{Zoom of the probability to explain the data set in Fig. \ref{speed} of the subsets from $l=6400$ to $l=6750$ of 14, 15, and 16 basis functions chosen in the $\{1, x, y, ...\, P_r(x) P_s(x), ... \}$ list, where $0 \le r+s \le 5$. The subsets are numbered with the shortest first and the later elements in the list omitted first. The probability of the remaining subsets is zero for all practical purposes.} \label{subsets}
\end{figure}

In order to carry out a more detailed analysis, regressions were carried out also by using the 190893 subsets of 14, 15, and 16 basis functions chosen in the $\{ 1, x, y, ...\, P_r(x) P_s(y), ... \}$ list, where $0 \le r+s \le 5$ and $P_r(x)$ is a Legendre polynomials of degree $r$. The results are shown in Fig.\ \ref{subsets}. The 14 basis functions whose linear combination -- which corresponds to a fifth degree polynomial -- most likely explains the data are $\{1$, $P_1(x)$, $P_1(y)$, $P_2(x)$, $P_2(y)$, $P_1(x)P_1(y)$, $P_3(x)$, $P_3(y)$, $P_1(x)P_2(y)$, $P_4(x)$, $P_4(y)$, $P_1(x)P_3(y)$, $P_2(x)P_2(y)$, $P_1(x)P_4(y)\}$.

A fallout of this Bayesian analysis are the probabilities of all the sets of smoothing polynomials considered to model the data. Consequently, when, as in this case, a number of basis function sets have a significant probability to explain the data, the quantities of interest -- in this case, the speed of sound values -- and the uncertainty inherent the model selection can be inferred by model averaging \cite{Wasserman:2000,Mana:2012}.

\section{Conclusions}
An analytical solution has been found for the problem of finding what basis functions must be used in linear regression analyses. It relies on the Bayesian model selection and complements the numerical results given in \cite{Gull:1988}. It uses the probability algebra to select among different basis function sets and encodes a preference for the smallest set capable to explain the data. In practice, a probability density is assigned to the regression coefficients prior the data are available, consistently with the given prior information and according to the maximum entropy principle. Next, the probability algebra allows this probability distribution to be updated according to the additional information delivered by the data. The regression probability is proportional to the normalising factor of the parameter likelihood times the parameter prior distribution.

A feature of this solution is that, for a large data sample, the regression probability depends only on the residuals and the number of free parameters. The smaller are the residuals, the higher the probability; but, a penalty exists for increasing the number of parameters. In addition, if a basis-function set explains the data exactly, its probability to explain the data is one.

\ack
This work was jointly funded by the European Metrology Research Programme (EMRP) participating countries within the European Association of National Metrology Institutes (EURAMET) and the European Union.

\end{document}